\documentclass[11pt]{article}
\usepackage{amsmath}
\usepackage{amssymb}
\usepackage{enumerate}
\usepackage{color}
\usepackage{xfrac}
\usepackage{bbm}
\usepackage{hyperref}
\usepackage{mathrsfs}
\usepackage{nicefrac}

\newtheorem{thm}{Theorem}[section]

\newtheorem{lemma}[thm]{Lemma}
\newtheorem{prop}[thm]{Proposition}
\newtheorem{rem}[thm]{Remark}
\newtheorem{cor}[thm]{Corollary}
\newtheorem{ex}[thm]{Example}

\newcommand{\bl}{\mathcal B_{\lambda}}
\newcommand{\tla}{\mathcal T_{\lambda,\alpha}}

\newcommand{\aaa}{\alpha}
\newcommand{\sss}{\sigma}

\begin{document}

\title{Approximations of generalized Bernstein functions}

\author{Stamatis Koumandos and
Henrik L. Pedersen\footnote{Research supported by grant DFF-1026-00267B from The Danish Council for Independent Research $|$ Natural Sciences}
}

\date{\today}
\maketitle

\begin{abstract}
We establish sharp inequalities involving the incomplete Beta and Gamma functions. These inequalities arise in the approximation of generalized Bernstein functions by  higher order Thorin-Bernstein functions. Furthermore, new properties of a related function, namely  $x^{\lambda}\Gamma(x)/\Gamma(x+\lambda)$ are derived.
\end{abstract}


\section{Introduction}

The results in this paper arise in the study of generalized Bernstein functions
and higher order Thorin-Bernstein functions. Let us begin by recalling the definitions.

Throughout this paper, $\lambda$ denotes a strictly positive number. A positive (meaning non-negative) function $f$ defined on $(0,\infty)$ having derivatives of all orders is a generalized Bernstein function of order $\lambda$ if $f'(x)x^{1-\lambda}$ is a completely monotonic function. This is equivalent to $f$ admitting an integral representation of the form
\begin{equation}\label{eq:int-rep-bl}
f(x)=ax^{\lambda}+b+\int_0^{\infty}\gamma(\lambda,xt)\frac{d\mu(t)}{t^{\lambda}},
\end{equation}
where $a$ and $b$ are non-negative numbers, and $\mu$ is a positive measure on $(0,\infty)$ making the integral converge for all $x>0$. Here, $\gamma(\lambda,x)$ denotes the incomplete gamma function defined by
$$
\gamma(\lambda,x)=\int_0^xe^{-t}t^{\lambda-1}\, dt.
$$
The class of generalized Bernstein functions of order $\lambda$ is denoted by $\bl$ and was studied in \cite{KP1}.  In view of  \eqref{eq:int-rep-bl}
 the incomplete gamma function can  be considered as the most important member of this family.
 In the case where $\lambda=1$ this is the class of ordinary Bernstein functions denoted by $\mathcal B$, see \cite{S}.

An important subclass of $\bl$ is the class  $\mathcal T_{\lambda,\alpha}$, of $(\lambda,\alpha)$-Thorin-Bernstein functions introduced and studied in \cite{KP4}. 
 A  function $f:(0,\infty)\to \mathbb R$ is in $\mathcal T_{\lambda,\alpha}$ if 
 \begin{equation}
 \label{eq:def-tb}
 f(x)=a x^{\lambda}+b+\int_0^{\infty}\gamma(\lambda,xt)\varphi(t)\, dt,
\end{equation}
where $a$ and $b$ are non-negative numbers, and $\varphi$ is a completely monotonic function of order $\alpha$.  A function $\varphi$ defined on $(0,\infty)$ is completely monotonic of order $\alpha$ if $x^{\alpha}\varphi(x)$ is completely monotonic.
For given $\lambda>0$, $\alpha$ is required to satisfy $\alpha<\lambda +1$. (Otherwise the integral in \eqref{eq:def-tb} will diverge unless $\varphi$ is identically equal to zero.)
The special case $\mathcal T_{1,1}$ is the class of Thorin-Bernstein functions  studied in \cite[Chapter 8]{S}.

In \cite[Proposition 2.2]{KP4} it was obtained that any function from the class $\bl$ can be approximated pointwise by a sequence of functions from $\cup_{\alpha<\lambda+1}\mathcal{T}_{\lambda,\alpha}$. 

One of the objectives of the present paper is to obtain a concrete version of this result. In Proposition \ref{dense} we show that $f\in \bl$ with representation \eqref{eq:int-rep-bl} is approximated pointwise by $\{f_n\}$, where $f_n\in\mathcal{T}_{\lambda,1-n}$ is given by 
$$
f_n(x)=ax^{\lambda}+b+\frac{\Gamma(\lambda+n)}{\Gamma(n)}\,\int_0^{\infty}  B\Big(\lambda, n, \frac{x}{x+n/t}\Big)\,\frac{d\mu(t)}{t^{\lambda}}.
$$
Here $B(a_1,a_2;x)$ is the incomplete Beta function, defined
 for $x\in [0,1)$, $a_1>0$, and $a_2\in \mathbb R$ as
$$
B(a_1,a_2;x)=\int_0^x t^{a_1-1}(1-t)^{a_2-1}\, dt.
$$

The proof of Proposition \ref{dense} requires applications of the dominated convergence theorem and these necessitated inequalities involving the incomplete Beta and Gamma functions. We believe that these inequalities are of independent interest and give them in the next result. 
\begin{prop}
 \label{thm:A1} $(\rm i)$  For  $\lambda<1$ and $x,y>0$, we have
\begin{equation}\label{eq:I}
\frac{\Gamma(\lambda+y)}{\Gamma(y)}\, B\Big(\lambda, y, \frac{x}{x+y}\Big)<\gamma(\lambda,x)\,.
\end{equation}
$(\rm ii)$  For  $\lambda>1$ and $x,y>0$, we have
  \begin{equation}\label{eq:II}
 y^{\lambda}\, B\Big(\lambda, y, \frac{x}{x+y}\Big)<\gamma(\lambda,x)\,.
\end{equation}
\end{prop}

It is convenient for us to introduce the function 
\begin{equation}
\label{eq:g-lambda}g_{\lambda}(x)=\frac{x^{\lambda}\Gamma(x)}{\Gamma(\lambda+x)}, \;\; x>0,
\end{equation}
which will play an important role in our investigations.

When $\lambda=1$,  \eqref{eq:I} coincides with \eqref{eq:II}, and an easy calculation shows that these inequalities reduce to the
elementary inequality
$
1+x/y<e^{x/y}$
for $x,y>0$.
We observe also that for any positive $\lambda$,
 \begin{equation*}
\lim_{x\to \infty}\Bigg[\frac{\Gamma(\lambda+y)}{\Gamma(y)}\, B\Big(\lambda, y, \frac{x}{x+y}\Big)\Bigg]=\lim_{x\to \infty}\gamma(\lambda,x)=\Gamma(\lambda)\,,
\end{equation*}
while inequality \eqref{eq:II}, as $x\to \infty$,  reduces to the sharp inequality $g_{\lambda}(y)<1$ for $y>0$ and $\lambda>1$,
see Corollary \ref{inq:gam}, below.

We note that both \eqref{eq:I} and \eqref{eq:II} are sharp in the sense that
\begin{equation}\label{eq:L}
\lim_{y\to \infty}\frac{\Gamma(\lambda+y)}{\Gamma(y)}\, B\Big(\lambda, y, \frac{x}{x+y}\Big)=\lim_{y\to \infty}y^{\lambda}\, B\Big(\lambda, y, \frac{x}{x+y}\Big)=\gamma(\lambda,x)\,.
\end{equation}
To see that \eqref{eq:L} holds, use the well-known fact $g_{\lambda}(y)\to 1$ as $y\to \infty$
(cf. \cite[6.1.47]{as}) and write
$$
B\left(\lambda,y; \frac{x}{x+y}\right)=y^{y}\,\int_0^{x}\frac{t^{\lambda-1}}{(y+t)^{\lambda+y}}\,dt.
$$
Since
$$
\lim_{y\to \infty}\Big(\frac{y}{y+t}\Big)^{\lambda+y}=e^{-t}\,,
$$
 an application of the dominated convergence theorem yields \eqref{eq:L}.

It turns out that the inequalities in Proposition \ref{thm:A1} are special cases of Theorem \ref{thm:A2}, and this result is one of  our main objectives in this work. Before stating it we give some notation.

Let $\{dP_{y,x}\}_{x,y>0}$ denote the family of gamma distributions with shape parameter $y$ and scale parameter $x/y$, meaning that $dP_{y,x}$ has density 
$$
\frac{1}{\Gamma(y)}\,(y/x)^y\,\, t^{y-1}\,e^{-ty/x}
$$
w.r.t.\ Lebesgue measure on the positive line. The distribution function  $P_{y,x}$ is given 
as $P_{y,x}(t)
=\gamma(y, yt/x)/\Gamma(y)$.
It is well-known and easy to see that 
$$
dP_{y,x}\rightarrow \epsilon_{x}, \quad \text{as}\quad y\to \infty,
$$
weakly, where $\epsilon_{x}$ is the point mass at $x$. Therefore, if $f\in \bl$  and $f$ is bounded, we have
$$
\lim_{y \to \infty}\int_{0}^{\infty} f(t)\,dP_{y,x}(t)=f(x)\;.
$$ 
The result can now be stated as follows.
\begin{thm}
 \label{thm:A2}
Let $f$ be any bounded and non-constant function in $\bl$.  
\begin{enumerate}
\item[(i)] 
For  $\lambda\leq1$ and $x,y>0$, we have
\begin{equation}\label{eq:III}
\int_{0}^{\infty} f(t)\,dP_{y,x}(t)<f(x)\;.
\end{equation}
\item[(ii)]For  $\lambda>1$ and $x,y>0$, we have
  \begin{equation}\label{eq:IV}
 g_{\lambda}(y)\,\int_{0}^{\infty} f(t)\,dP_{y,x}(t)< f(x)\;,
\end{equation}
where $g_{\lambda}$ is defined in \eqref{eq:g-lambda}.
\end{enumerate}
\end{thm}

The last part of the paper is devoted to a more detailed study of the function $g_{\lambda}$ in \eqref{eq:g-lambda} for any $\lambda >0$. See Proposition \ref{prop:extra-g-lambda}, Proposition \ref{prop:logg-lambda-in-s2} and Proposition \ref{prop:-logg-lambda-in-s2}.

\section{Main results}
The proof of Theorem \ref{thm:A2} relies on the next Lemmas.
\begin{lemma}
 \label{lemma:P1}
 For $x,y>0$ the following formula holds true
 \begin{equation}\label{eq:BG1}
 \frac{\Gamma(\lambda+y)}{\Gamma(y)}\, B\Big(\lambda, y, \frac{x}{x+y}\Big)=\Gamma(\lambda)-x^{\lambda}\, \int_{0}^{\infty} \frac{\gamma(y,\, yt)}{\Gamma(y)}\,
 t^{\lambda-1}e^{-xt}\, dt\,.
 \end{equation}
 \end{lemma}
{\it Proof.} We know from \cite[Proposition  2.8]{KP4} that
 for any positive measure $\mu$,$\;$  $\alpha<\lambda+1$ and $x>0$
\begin{align*}
\int_0^{\infty}\gamma(\lambda,xt)t^{-\alpha}\mathcal L (\mu)(t)\, dt=\Gamma(\lambda+1-\alpha)\int_0^{\infty}B\left(\lambda,1-\alpha; \frac{x}{x+s}\right)\, \frac{d\mu(s)}{s^{1-\alpha}}\,,
\end{align*}
where $\mathcal L (\mu)$
is the Laplace transform of $\mu$. Taking $\alpha=1-y$ and $\mu=\epsilon_{y}$ the formula becomes
\begin{equation}\label{eq:BG2}
 \frac{\Gamma(\lambda+y)}{\Gamma(y)}\, B\Big(\lambda, y, \frac{x}{x+y}\Big)=\frac{y^y}{\Gamma(y)}\,\int_0^{\infty}\gamma(\lambda,xt)\,t^{y-1}\, e^{-yt}\, dt\,.
 \end{equation}
Integrating by parts on the right hand side of \eqref{eq:BG2} we obtain \eqref{eq:BG1}. \hfill $\square$

\begin{lemma}
 \label{lemma:P2}
 For $y>0$ the following formula holds true
 \begin{equation}\label{eq:GG1}
 \int_{0}^{1} \frac{\gamma(y,\, yt)}{\Gamma(y)}\,t^{\lambda-1}\,dt+ \int_{1}^{\infty} \left(\frac{\gamma(y,\, yt)}{\Gamma(y)}-1\right)\,t^{\lambda-1}\,dt=
 \frac{1}{\lambda}\,\left(1-\frac{1}{g_{\lambda}(y)}\right)\,.
 \end{equation}

 \end{lemma}
{\it Proof.} Integrating by parts we obtain
$$
 \int_{0}^{1} \frac{\gamma(y,\, yt)}{\Gamma(y)}\,t^{\lambda-1}\,dt= \frac{1}{\lambda}\,\frac{\gamma(y,\, y)}{\Gamma(y)}-\frac{y^y}{\lambda\,\Gamma(y)}\, \int_{0}^{1}t^{\lambda+y-1}\,  e^{-yt}\, dt\,.
$$
It is elementary to show that 
$$
\lim_{t\to \infty}\left(\frac{\gamma(y,\, yt)}{\Gamma(y)}-1\right)\,t^{\lambda}=0\,.
$$
Therefore, an integration by parts gives
$$
 \int_{1}^{\infty} \left(\frac{\gamma(y,\, yt)}{\Gamma(y)}-1\right)\,t^{\lambda-1}\,dt=\frac{1}{\lambda}-\frac{1}{\lambda}\,\frac{\gamma(y,\, y)}{\Gamma(y)}-\frac{y^y}{\lambda\,\Gamma(y)}\, \int_{1}^{\infty}t^{\lambda+y-1}\,  e^{-yt}\, dt\,.
$$
Adding the above integrals we arrive at
\begin{align*}
&\int_{0}^{1} \frac{\gamma(y,\, yt)}{\Gamma(y)}\,t^{\lambda-1}\,dt+ \int_{1}^{\infty} \left(\frac{\gamma(y,\, yt)}{\Gamma(y)}-1\right)\,t^{\lambda-1}\,dt\\[+5pt]
&=\frac{1}{\lambda}-\frac{y^y}{\lambda\,\Gamma(y)}\, \int_{0}^{\infty}t^{\lambda+y-1}\,  e^{-yt}\, dt=
\frac{1}{\lambda}\,\left(1-\frac{\Gamma(\lambda+y)}{y^{\lambda}\,\Gamma(y)}\right)=\frac{1}{\lambda}\left(1-\frac{1}{g_{\lambda}(y)}\right)\,.
\end{align*}
The proof is complete. \hfill $\square$

The following monotonicity properties of the function $g_{\lambda}$ in \eqref{eq:g-lambda} are well known,  see \cite{ism} and \cite{KL}.
\begin{lemma}
 \label{lemma:P3}
 For  $\lambda<1$, $g_{\lambda}$ is strictly decreasing on $(0, \, \infty)$, and for  $\lambda>1$, $g_{\lambda}$ is strictly increasing on $(0, \, \infty)$.
 \end{lemma}
 In view of this lemma and the previously mentioned fact 
 $\lim_{y\to \infty}g_{\lambda}(y)=1$ we immediately deduce the following.
\begin{cor}\label{inq:gam}
For $\lambda<1$, we have $g_{\lambda}(x)>1$ for $x>0$; for
$\lambda>1$, we have $g_{\lambda}(x)<1$ for $x>0$.
 \end{cor}
 As mentioned in the introduction the results on $g_{\lambda}$ are strengthened in last part of Section \ref{sec:app}.
 \medskip
 
\noindent
{\it Proof of Theorem \ref{thm:A2}.}  Let $f(\infty):=\lim_{t\to \infty}f(t)$ ($f$ being increasing and bounded). Furthermore, denoting by $\chi_{[x,\infty)}$ the indicator function of the interval $[x,\infty)$,  
we may write 
\begin{equation*}
f(x)=f(\infty)-\int_{0}^{\infty} \chi_{[x,\infty)}(t)\,f^{\prime}(t)\, dt\,.
\end{equation*}
Integrating by parts we obtain
$$
\int_{0}^{\infty} f(t)\,dP_{y,x}(t)=f(\infty)-\int_{0}^{\infty} P_{y,x}(t)\,f^{\prime}(t)\, dt\,.
$$
Hence, \eqref{eq:III} is equivalent to
$$
\int_{0}^{\infty}\big( P_{y,x}(t)-\chi_{[x,\infty)}(t)\big)\,f^{\prime}(t)\, dt>0.
$$

Taking into consideration that $t^{1-\lambda}\,f^{\prime}(t)$ is completely monotonic and   using Lemma \ref{lemma:P2} and Lemma \ref{lemma:P3},  we obtain
\begin{align*}
&\int_{0}^{\infty}\big( P_{y,x}(t)-\chi_{[x,\infty)}(t)\big)\,f^{\prime}(t)\, dt\\[+5pt]
&=\int_{0}^{x} P_{y,x}(t)\,f^{\prime}(t)\, dt+\int_{x}^{\infty}\big(P_{y,x}(t)-1\big)\,f^{\prime}(t)\, dt\\[+5pt]
&\geq x^{1-\lambda}\,f^{\prime}(x)\,\left(\int_{0}^{x} P_{y,x}(t)\,t^{\lambda-1}\, dt+\int_{x}^{\infty}\big( P_{y,x}(t)-1\big)\,t^{\lambda-1}\, dt\right)\\[+5pt]
&= x^{1-\lambda}\,f^{\prime}(x) \, x^{\lambda}  \, \frac{1}{\lambda}\,\left(1-\frac{1}{g_{\lambda}(y)}\right).
\end{align*}
This is positive, by Corollary \ref{inq:gam}, and \eqref{eq:III} is thus established. 

For the proof of \eqref{eq:IV} we observe that this is equivalent to
\begin{equation}\label{eq:YX}
f(\infty)\,\big(1-g_{\lambda}(y)\big)+\int_{0}^{\infty} \big(g_{\lambda}(y)P_{y,x}(t)- \chi_{[x,\infty)}(t)\big)\,f^{\prime}(t)\, dt>0\,.
\end{equation}
We write
\begin{align*}
&\int_{0}^{\infty} \big(g_{\lambda}(y)P_{y,x}(t)-\chi_{[x,\infty)}(t)\big)\,f^{\prime}(t)\, dt\\[+5pt]
&=g_{\lambda}(y)\,\left(\int_{0}^{x} P_{y,x}(t)\,f^{\prime}(t)\, dt+\int_{x}^{\infty}\big(P_{y,x}(t)-1\big)\,f^{\prime}(t)\, dt\right)\\[+5pt]
&+\big(g_{\lambda}(y)-1\big)\,\big(f(\infty)-f(x)\big)\,.
\end{align*}
As in the previous case, we have
\begin{align*}
\int_{0}^{x} P_{y,x}(t)\,f^{\prime}(t)\, dt+\int_{x}^{\infty}\big(P_{y,x}(t)-1\big)\,f^{\prime}(t)\, dt \geq \frac{x}{\lambda}\,f^{\prime}(x) \, 
\left(1-\frac{1}{g_{\lambda}(y)}\right)\,.
\end{align*}
Combining the above, the left hand side of \eqref{eq:YX} becomes
\begin{align*}
&f(x)\,\big(1-g_{\lambda}(y)\big)+g_{\lambda}(y)\,\left(\int_{0}^{x} P_{y,x}(t)\,f^{\prime}(t)\, dt+\int_{x}^{\infty}\big(P_{y,x}(t)-1\big)\,f^{\prime}(t)\, dt\right)\\[+5pt]
& \geq f(x)\,\big(1-g_{\lambda}(y)\big)+\frac{x}{\lambda}\,f^{\prime}(x) \, \big(g_{\lambda}(y)-1\big)\\[+3pt]
&= \big(1-g_{\lambda}(y)\big)\,\left(f(x)-\frac{x\,f^{\prime}(x)}{\lambda}\right)>0.
\end{align*}
We have used again Lemma \ref{lemma:P3}, Corollary \ref{inq:gam} and the fact that  for any function $f\in \bl$, $t^{-\lambda}f(t)$ is completely monotonic (see \cite[[Corollary 2.1]{KP1}), and in fact strictly decreasing when $f$ is bounded.

The proof of the theorem is complete. \hfill $\square$

\begin{rem}
    The inequality in Theorem \ref{thm:A2} when $\lambda\leq1$ actually holds for all bounded, positive and increasing functions $f$ for which $t^{1-\lambda}f'(t)$ is strictly decreasing. For $\lambda>1$ we also need that $t^{-\lambda}f(t)$ is strictly decreasing. 
    
\end{rem}
\noindent
{\it Proof of Proposition \ref{thm:A1}.} Taking $x=1$ and $f(t)=\gamma(\lambda,zt)\in \bl$, $z>0$ in Theorem \ref{thm:A2},  and using \eqref{eq:BG2} we obtain the assertions.\hfill $\square$ 

\medskip


\begin{cor}
\label{cor:cm}
Let $f\in \bl$ be bounded and non-constant, and let $\omega$ be the measure in the Bernstein representation of the completely monotonic function $t^{-\lambda}\,f(t)$. Then  
\begin{enumerate}  
\item [(i)]  
For  $\lambda\leq1$ and $x,y>0$, we have
\begin{equation*}
y^{y}\, \frac{\Gamma(\lambda+y)}{\Gamma(y)}\,\int_{0}^{\infty} \frac{d\omega(s)}{(y+xs)^{\lambda+y}}<\mathcal L (\omega)(x)\;.
\end{equation*}
\item[(ii)]
For  $\lambda>1$ and $x,y>0$, we have
  \begin{equation*}
 y^{\lambda+y}\, \,\int_{0}^{\infty} \frac{d\omega(s)}{(y+xs)^{\lambda+y}}<\mathcal L (\omega)(x)\;.
\end{equation*}
\end{enumerate}
\end{cor} 
Here, $\mathcal L$ stands for the Laplace transform.

{\it Proof.}
This is immediate from \cite[Corollary 2.6]{KP4} and Theorem \ref{thm:A2}. \hfill $\square$
\medskip

It is well-known (see e.g.\ \cite[Proposition 3.12]{S}) that $f\in \mathcal {B}_1$ is bounded if and only if there is $c\geq 0$ and a bounded completely monotonic function $g$ such that $f=c-g$. Theorem \ref{thm:A2} thus yields the following corollary.
\begin{cor}
\label{cor:bdd_cm}
    Let $g$ be any bounded and non-constant completely monotonic function.  
 For $x,y>0$, we have
$$
g(x)<\int_{0}^{\infty} g(t)\,dP_{y,x}(t)=\frac{(y/x)^y}{\Gamma(y)}\int_0^{\infty}g(t)t^{y-1}e^{-ty/x}\, dt.
$$
\end{cor}

As mentioned in the introduction the motivation for our results is a concrete approximation result, stated in Proposition \ref{dense}.
Before proving it we recall for the reader's convenience a characterization 
in \cite[Corollary 3.2]{KP4}:
 A function $f$ belongs to  $\tla$ if and only if  $f$ can be represented as
  $$
 f(x)=ax^{\lambda}+b+\Gamma(\lambda+1-\alpha)\int_0^{\infty}B\left(\lambda,1-\alpha;\frac{x}{x+t}\right)\frac{d\mu(t)}{t^{1-\alpha}},
 $$
 where $a$, $b$ are non-negative numbers and $\mu$ is a positive measure making the integral converge.
\begin{prop}\label{dense}
Let $f\in \bl$ have the representation \eqref{eq:int-rep-bl}. The function 
$$
f_n(x)=ax^{\lambda}+b+\frac{\Gamma(\lambda+n)}{\Gamma(n)}\,\int_0^{\infty}  B\Big(\lambda, n, \frac{x}{x+n/t}\Big)\,\frac{d\mu(t)}{t^{\lambda}}
$$
belongs to $\mathcal T_{\lambda,1-n}$ and 
$$
\lim_{n\to \infty} f_{n}(x)=f(x), \quad \text{for all} \quad x>0.
$$
\end{prop}
{\it Proof.} Suppose that $f$ has the representation \eqref{eq:int-rep-bl}, let $\kappa(t)=n/t$ and let $\omega=\kappa(t^{-\lambda-n}d\mu(t))$ denote the image measure of $t^{-\lambda-n}d\mu(t)$. Since
$$
\int_0^{\infty}B\left(\lambda, n, \frac{x}{x+n/t}\right)\,\frac{d\mu(t)}{t^{\lambda}}=
n^n\int_0^{\infty} B\left(\lambda, n, \frac{x}{x+s}\right)\,\frac{d\omega(s)}{s^{n}}
$$
we see that $f_n$ given in the statement of the proposition belongs to $\mathcal{T}_{\lambda,1-n}$.

For $\lambda<1$ we get from \eqref{eq:I}, \eqref{eq:L} and the dominated convergence theorem, that $f_n(x)\to f(x)$ as $n\to \infty$.
In the case where $\lambda>1$,  we argue similarly using \eqref{eq:II} instead of \eqref{eq:I}.
When $\lambda=1$ we use the dominated convergence theorem and the elementary inequality
$(1+x/n)^{n}<e^{x}$ for $x>0$.

The proof is complete. \hfill $\square$
\section{On the auxiliary function}
\label{sec:app}





This section is devoted to some additional properties of the function $g_{\lambda}$ defined in \eqref{eq:g-lambda}. In this investigation the class $\mathcal S_{\lambda}$ of generalized Stieltjes functions of order $\lambda$ plays an important role: $g\in \mathcal S_{\lambda}$ if 
$$
g(x)=\int_0^{\infty}\frac{d\mu(t)}{(t+x)^{\lambda}}+c
$$
for a positive measure $\mu$ (making the integral converge) and a non-negative constant $c$.

There is a connection between these functions and the higher order Thorin-Bernstein functions: $f\in \mathcal{T}_{\lambda,\alpha}$ if and only if $x^{\lambda-1}f'(x)\in \mathcal S_{\lambda+1-\alpha}$. (See \cite[Theorem 2.14]{KP4}.)

In Proposition \ref{prop:glambda-t} it is established that $g_{\lambda}\in \mathcal{T}_{\lambda-1,-3}$ when $\lambda>1$. Let us begin by giving a weaker, preliminary result, which, however, is an extension of Lemma \ref{lemma:P3}.
 \begin{prop}
 \label{prop:extra-g-lambda}
 For $\lambda>1$ we have $g_{\lambda}\in \mathcal{B}_{\lambda-1}\setminus \cup_{\tau<\lambda-1}\mathcal{B}_{\tau}$.
 \end{prop}
{\it Proof.} Differentiating we get
$$
g_{\lambda}^{\prime}(x)=g_{\lambda}(x)\,\left(\,\frac{\lambda}{x}+\psi(x)-\psi(x+\lambda)\right),
$$
where $\psi(x) =
\Gamma'(x)/\Gamma(x)$ is the psi or digamma function. The function $\sigma_{\lambda}(x):=\lambda/x+\psi(x)-\psi(x+\lambda)$ is seen to have the representation
\begin{equation*}
\sigma_{\lambda}(x)=\int_{0}^{\infty} e^{-xt}\, \Big(\lambda-\frac{1-e^{-\lambda t}}{1-e^{-t}}\Big)\,dt.
\end{equation*}
Now $t\mapsto \lambda-(1-e^{-\lambda t})/(1-e^{-t})$ is positive and increasing on $[0,\infty)$, and takes the value zero at the origin.  Hence not only $\sigma_{\lambda}(x)$ but also  $x\sigma_{\lambda}(x)$ is completely monotonic (see e.g.\ \cite[Theorem 1.3]{KP0}) and
\begin{equation}
\label{eq:xsigma-laplace}
x\sigma_{\lambda}(x)=\mathcal L\left(\left(\lambda-\frac{1-e^{-\lambda t}}{1-e^{-t}}\right)'\right)(x).
\end{equation}
It is well-known and easy to see (cf.\ \cite{ism}) that $\Gamma(x+1)/\Gamma(x+\lambda)$ is a completely monotonic function and since 
\begin{equation}\label{eq:der}
x^{2-\lambda}\,g_{\lambda}^{\prime}(x)=x\sigma_{\lambda}(x)\frac{\Gamma(x+1)}{\Gamma(x+\lambda)},
\end{equation}
it follows that $x^{2-\lambda}\,g_{\lambda}^{\prime}(x)$  is completely monotonic, i.e.\
that $g_{\lambda}\in \mathcal{B}_{\lambda-1}$. 

If $g_{\lambda}\in\mathcal{B}_{\tau}$ for some $\tau<\lambda-1$  then $g_{\lambda}(x)/x^{\tau}$ would be completely monotonic (cf.\ \cite[Corollary 3.8, (iv)]{S} or \cite[Corollary 2.1]{KP1}). 
Since by inspection $g_{\lambda}(x)/x^{\tau}\to 0$ as $x\to 0^+$ this would force $g_{\lambda}$ to be the zero function, which is a contradiction. 

The proof of the proposition is complete. \hfill $\square$

\medskip

Before a more detailed treatment in the situation where $\lambda>1$, let us investigate the case $\lambda\in (0,1)$. 
Here it is known that (see \cite{ism}) $g_{\lambda}$ is logarithmically completely monotonic. We strengthen and complement this result in Proposition \ref{prop:logg-lambda-in-s2}.

For the reader's convenience we recall that $f:(0,\infty)\to (0,\infty)$ is logarithmically completely monotonic if $-f'(x)/f(x)$ is completely monotonic.

Any function from $\mathcal S_2$ is logarithmically completely monotonic. This is  a deep result of Kristiansen (\cite{kristiansen}), see also \cite[Theorem 2.1]{BKP}. Proposition \ref{prop:logg-lambda-in-s2} thus entails that when $\lambda<1$, $g_{\lambda}$ and $\log g_{\lambda}$ are logarithmically completely monotonic functions. 
\begin{prop}
\label{prop:logg-lambda-in-s2}
    For $\lambda\in (0,1)$ the function $\log g_{\lambda}$ belongs to $\mathcal S_2$ but not to any of the Stieltjes classes $\mathcal S_{\tau}$ for $\tau<2$.
\end{prop}

    \noindent
{\it Proof of Proposition \ref{prop:logg-lambda-in-s2}.}
From \cite[(6.10)]{KL} we have 
    $$
    -(\log g_{\lambda})'(x)=\int_0^{\infty}e^{-xu}\Phi(u)\, du,
    $$
    where 
    $$
    \Phi(u)=\frac{1-e^{-\lambda u}}{1-e^{-u}}-\lambda.
    $$
    Using Fubini's theorem and the relation $\lim_{u\to \infty}g_{\lambda}(u)=1$ we get
    $$
    \log g_{\lambda}(x)=\int_0^{\infty}e^{-xu}u\frac{\Phi(u)}{u^2}\, du.
    $$
    Thus $\log g_{\lambda}$ belongs to $\mathcal S_2$ if and only if $\Phi(u)/u^2$ is completely monotonic. Now, $\Phi=\mathcal L(\nu)$, where
    $$
    \nu=\sum_{k=0}^{\infty}(\epsilon_k-\epsilon_{\lambda+k})-\lambda\epsilon_0.
    $$
    This gives 
    $$
    \frac{\Phi(u)}{u^2}=\mathcal L(\xi d\xi\ast \nu)(u),
    $$
    and it is easily seen that $\xi d\xi\ast \nu$ has density 
    $$
    \varphi(\xi)=\sum_{k=0}^{\infty}\left(\chi_{[k,\infty)}(\xi)(\xi-k)-\chi_{[k+\lambda,\infty)}(\xi)(\xi-k-\lambda)\right) - \lambda \xi
    $$
    w.r.t\ Lebesgue measure on $[0,\infty)$. We claim that $\varphi$ is non-negative on $[0,\infty)$. Indeed a computation shows that $\varphi$ is $1$-periodic on $[0,\infty)$ and for $\xi \in [0,\lambda)$ it follows by inspection that
    $\varphi(\xi)=(1-\lambda)\xi\geq 0$. In the case where $\xi\in [\lambda,1)$ we similarly find $\varphi(\xi)=\lambda(1-\xi)\geq 0$.    
    The complete monotonicity of $\Phi(u)/u^2$ is verified, and hence $\log g_{\lambda}\in \mathcal S_{2}$.
    
    In the same way, it follows that $\log g_{\lambda}\in \mathcal S_{1+\sigma}$ if and only if 
$$
u\mapsto \frac{1}{u^{1+\sigma}}\left(\lambda-\frac{1-e^{-\lambda u}}{1-e^{-u}}\right)
$$
is completely monotonic. This is equivalent to the function $\varphi_{\sigma}$ given by
$$
\varphi_{\sigma}(\xi)=
\sum_{k=0}^{\infty}\left(\chi_{[k,\infty)}(\xi)(\xi-k)^{\sigma}-\chi_{[k+\lambda,\infty)}(\xi)(\xi-k-\lambda)^{\sigma}\right)-\lambda \xi^{\sigma}
$$
being non-negative on $[0,\infty)$. However, computations show that $\varphi_{\sigma}(1)<0$ and $\varphi_{\sigma}(\lambda)>0$ when $\sigma\in (0,1)$.\hfill $\square$
\medskip

Let us now turn to the case where $\lambda>1$.

For $\lambda\in (1,2]$ we have $g_{\lambda}\in \mathcal B$ by Propostion \ref{prop:extra-g-lambda}. Therefore, both $g_{\lambda}'$ and $1/g_{\lambda}$ are completely monotonic. Hence also the product is completely monotonic and this means that $1/g_{\lambda}$ is logarithmically completely monotonic. 

It is proved below in Proposition \ref{prop:-logg-lambda-in-s2} that $\log (1/g_{\lambda})$ belongs to $\mathcal S_2$ for $\lambda>1$. As noticed earlier this implies that $1/g_{\lambda}$ and its logarithm are also logarithmically completely monotonic.

\begin{prop}
\label{prop:-logg-lambda-in-s2}
    For $\lambda >1$ the function $-\log g_{\lambda}$ belongs to $\mathcal S_2$. When $\lambda>1$ is not an integer, $-\log g_{\lambda}$ does not belong to any $\mathcal S_{\tau}$ where $\tau<2$.
\end{prop}
(For the elementary integer case, see Remark \ref{rem:integer}.)

\noindent
{\it Proof.} Let $\sigma\in (0,1]$. It follows in the same way as in the proof of Proposition  \ref{prop:logg-lambda-in-s2} that $-\log g_{\lambda}\in \mathcal S_{1+\sigma}$ if and only if the function
$$
\varphi(\xi)=
\lambda\xi^{\sigma}-\sum_{k=0}^{\infty}\left(\chi_{[k,\infty)}(\xi)(\xi-k)^{\sigma}-\chi_{[k+\lambda,\infty)}(\xi)(\xi-k-\lambda)^{\sigma}\right)
$$
is non-negative on $[0,\infty)$. For $\sigma=1$ the non-negativity of $\varphi$ is verifed as above, so $-\log g_{\lambda}\in \mathcal S_2$. 

Now let $\sss\in (0,1)$. We write $\lambda=N+\alpha$, where $N\in \{1,2,\ldots\}$ and $\alpha\in (0,1)$. The idea is show that, as $n$ tends to infinity, $\varphi(n+\alpha)$ converges to a negative number.

For $n\geq N$ we have $\varphi(n+\aaa)=S_n(N)$, where
$$
S_n(N)=(N+\aaa)(n+\aaa)^{\sss}-\sum_{k=0}^n(k+\aaa)^{\sss}+\sum_{k=0}^{n-N}k^{\sss}.
$$
First of all notice that it is sufficient to consider $N=0$. Indeed
$$
S_n(N)-S_n(0)=N(n+\aaa)^{\sss}-\sum_{k=n-N+1}^nk^{\sigma}=
\sum_{k=n-N+1}^n\left((n+\aaa)^{\sss}-k^{\sigma}\right).
$$
Here, 
\begin{align*}
  \sum_{k=n-N+1}^n\left((n+\aaa)^{\sss}-k^{\sigma}\right)&=\sum_{k=n-N+1}^n\sss \int_k^{n+\aaa}t^{\sss-1}\, dt\\
  &\leq N\sss (n-N+1)^{\sss-1}(n+\aaa-n+N-1)\\
  &=N\sss (n-N+1)^{\sss-1}(N+\aaa-1)\\
  &=o(1).
\end{align*}
Hence $S_n(N)-S_n(0)=o(1)$ and it is enough to consider $S_n(0)$. Notice that also 
$(n+\aaa)^{\sss}-n^{\sss}=o(1)$ (a particular case of the relation above for $N=1$), so we have 
$$
S_n(0)=\aaa(n+\aaa)^{\sss}-\sum_{k=0}^{n-1}(k+\aaa)^{\sss}+\sum_{k=0}^{n-1}k^{\sss}+o(1).
$$
The difficulty is that each of the three terms in this expression tend to infinity as $n$ tends to infinity. To overcome this we apply the Euler-MacLaurin summation formula (see e.g.\ \cite[(D.1.3)]{aar}), 
$$
\sum_{k=m+1}^{n-1}f(k)=-\nicefrac{1}{2}(f(n)+f(m))+\int_m^nf(x)\, dx+\int_m^n(x-[x]-\nicefrac{1}{2})f'(x)\, dx
$$
to $f(x)=(x+\aaa)^{\sss}$ and $m=0$ and to $f(x)=x^{\sss}$ and $m=1$, giving
\begin{align*}
    S_n(0)&=o(1)+\left\{\aaa(n+\aaa)^{\sss}-\frac{1}{\sss+1}\left((n+\aaa)^{\sss+1}-n^{\sss+1}\right)\right\}\\
    &\phantom{=}+\frac{1}{2}\left((n+\aaa)^{\sss}-n^{\sss}\right)+\frac{1}{\sss+1}\left(\aaa^{\sss+1}-1\right)-\frac{1}{2}\left(\aaa^{\sss}-1\right)\\
    &\phantom{=}-\int_0^n(x-[x]-\nicefrac{1}{2})\sss (x+\aaa)^{\sss-1}\, dx+\int_1^n(x-[x]-\nicefrac{1}{2})\sss x^{\sss-1}\, dx.
\end{align*}
The expression in the curly brackets is $o(1)$, since
\begin{align*}
  \lefteqn{\aaa(n+\aaa)^{\sss}-\frac{1}{\sss+1}\left((n+\aaa)^{\sss+1}-n^{\sss+1}\right) }\\
  &=\aaa(n+\aaa)^{\sss}-\frac{1}{\sss+1}n^{\sss+1}\left(1+(\sss+1)\aaa/n+O(\nicefrac{1}{n^2})-1\right)\\
  &=\aaa(n+\aaa)^{\sss}-\aaa n^{\sss}+o(1)=o(1).
\end{align*}
The limit of each of the two integrals above, as $n$ tends to infinity, exists; this can be seen e.g.\ by partial integration. We thus have
\begin{align*}
\lim_{n\to \infty}S_n(0)&=
\frac{1}{\sss+1}\left(\aaa^{\sss+1}-1\right)-\frac{1}{2}\left(\aaa^{\sss}-1\right)\\
&\phantom{=}+\int_o^{\aaa}(t-\aaa-[t-\aaa]-\nicefrac{1}{2})\sss t^{\sss-1}\, dt-\int_0^1(t-[t]-\nicefrac{1}{2})\sss t^{\sss-1}\, dt\\
    &\phantom{=}-\int_{-\aaa}^{\infty}(x-[x]-\nicefrac{1}{2})\sss (x+\aaa)^{\sss-1}\, dx\\
    &\phantom{=}+\int_{-1}^{\infty}(x-[x]-\nicefrac{1}{2})\sss (x+1)^{\sss-1}\, dx.
\end{align*}
The two integrals in the second line are easily computed and cancel with right hand side members in the first line. This gives
\begin{align*}
\lim_{n\to \infty}S_n(0)
    &=-\int_{-\aaa}^{\infty}(x-[x]-\nicefrac{1}{2})\sss (x+\aaa)^{\sss-1}\, dx\\
    &\phantom{=}
    +\int_{-1}^{\infty}(x-[x]-\nicefrac{1}{2})\sss (x+1)^{\sss-1}\, dx.
\end{align*}
We claim that this expression is non-positive. (It is, by the way, equal to $\zeta(-\sss,1)-\zeta(-\sss,\aaa)$, where $\zeta$ denotes the Hurwitz zeta function, see e.g.\ \cite{dlmf}.)

By linear variable transformations it follows that
\begin{align*}
    \lim_{n\to \infty}S_n(0)&=
    -\sss \sum_{k=0}^{\infty}\int_0^1\left(1-\aaa+[t-1]-[t-\aaa]\right)(t+k)^{\sss-1}\, dt.
\end{align*}
Here, 
\begin{align*}
    \lefteqn{\int_0^1\left(1-\aaa+[t-1]-[t-\aaa]\right)(t+k)^{\sss-1}\, dt}\\
    &=(1-\aaa)\int_0^{\aaa}(t+k)^{\sss-1}\,dt-\aaa\int_{\aaa}^1(t+k)^{\sss-1}\,dt\\
    &>(1-\aaa)\aaa (\aaa+k)^{\sss-1}-\aaa (1-\aaa)(\aaa+k)^{\sss-1}=0.
\end{align*}
This shows finally, that for $\aaa,\sss \in (0,1)$, we have $\lim_{n\to \infty}S_n(0)<0$. \hfill $\square$

We complete this section by showing that $g_{\lambda}$ is a higher order Thorin-Bernstein function:
\begin{prop}
\label{prop:glambda-t}
 For $\lambda>1$, $g_{\lambda}\in \mathcal T_{\lambda-1,-3}$.  
\end{prop}
{\it Proof.} As explained in the beginning of this section, it is equivalent to prove that 
$x^{2-\lambda}g_{\lambda}'(x)\in \mathcal S_{\lambda+3}$. Since in fact $\Gamma(x+1)/\Gamma(x+\lambda)$ belongs to $\mathcal S_{\lambda}$, see \cite[Theorem 4.6]{BKP}, \eqref{eq:der} and Lemma \ref{lemma:xsigma} below yield that 
$$
x^{2-\lambda}g_{\lambda}'(x)\in \mathcal S_{\lambda}\cdot \mathcal S_3\subseteq \mathcal S_{\lambda+3}.
$$
This completes the proof.\hfill $\square$
\begin{lemma}
\label{lemma:xsigma}
    For $\lambda >1$, the function $x\sigma_{\lambda}(x)=\lambda+x(\psi(x)-\psi(x+\lambda))$ belongs to $\mathcal S_3$.
\end{lemma}
{\it Proof.} The relation \eqref{eq:xsigma-laplace} shows that it is sufficient to verify that the function $\Xi$ defined as 
$$
    \Xi(t)=-\frac{1}{t^2}\left(\frac{1-e^{-\lambda t}}{1-e^{-t}}\right)'
    $$
is completely monotonic. 

Since 
$$
-\left(\frac{1-e^{-\lambda t}}{1-e^{-t}}\right)'= \mathcal L\left(\sum_{k=0}^{\infty}(k\epsilon_k-(k+\lambda)\epsilon_{k+\lambda})\right)(t)
$$
and $1/t^2=\mathcal L(sds)(t)$ we get 
$$
\Xi(t)=\mathcal L\left(sds\ast \sum_{k=0}^{\infty}(k\epsilon_k-(k+\lambda)\epsilon_{k+\lambda})\right)(t)=\mathcal L(\phi)(t),
$$
where
$$
\phi(s)=\sum_{k=0}^{\infty}(k \chi_{[k,\infty)}(s)(s-k)-(k+\lambda)\chi_{[k+\lambda,\infty)}(s)(s-k-\lambda)).
$$
We prove that $\phi(s)\geq 0$ for $s\geq 0$ along the same lines as Proposition \ref{prop:logg-lambda-in-s2}. The recurrence relations 
\begin{enumerate}[(i)]
    \item $\phi(s+1)=\phi(s)+\phi_1(s)$, where $\phi_1(s)=\sum_{k=0}^{\infty}( \chi_{[k,\infty)}(s)(s-k)-\chi_{[k+\lambda,\infty)}(s)(s-k-\lambda))-\lambda \chi_{[\lambda-1,\infty)}(s)(s+1-\lambda)$,
\item  $\phi_1(s+1)=\phi_1(s)+\phi_2(s)$, where $\phi_2(s)=s+1+(\lambda-1) \chi_{[\lambda-1,\infty)}(s)(s+1-\lambda)-\chi_{[\lambda-2,\infty)}(s)(s+2-\lambda)$
\end{enumerate}
are easily seen to hold. By inspection, $\phi_2(s)\geq 0$ for $s\geq 0$, and also $\phi_1(s)\geq 0$ for $s\in[0,1]$. Thus $\phi_1(s)\geq 0$ for $s\geq 0$. Since $\phi(s)=0$ for $s\in [0,1]$ we see that also $\phi(s)\geq 0$ for $s\geq 0$.\hfill $\square$

\begin{rem}
\label{rem:integer}
    When $\lambda=N\in\{2,3,\ldots\}$ we have $g_N(x)=\prod_{k=1}^{N-1}\frac{x}{x+k}$ so that 
    \begin{align*}
     -\log g_N(x)&=\sum_{k=1}^{N-1}\log (1+k/x)\in \mathcal S_1,\\
     \sigma_N(x)&=\sum_{k=1}^{N-1}\frac{k}{x(x+k)}\in \mathcal S_2,\\ 
     x\sigma_N(x)&=\sum_{k=1}^{N-1}\frac{k}{x+k}\in \mathcal S_1,\ \text{and}\\
     x^{2-N}g_N'(x)&=\left(\prod_{k=1}^{N-1}\frac{1}{x+k}\right)\sum_{j=1}^{N-1}\frac{j}{j+x}\in  \mathcal S_N.
    \end{align*} 
    In particular, $g_N\in  \mathcal T_{N-1,0}$. 
\end{rem}
\section{Asymptotic expansions of bounded  $\mathcal T_{\lambda,0}$ functions}
\label{sec:tbd}

Let $f$ be a bounded complete generalized Bernstein function of order $\lambda$. According to \cite[Proposition 2.2 and (11)]{KP1}, $f$ has the representation
\begin{equation}
    \label{eq:bounded}
f(x)=b+\int_0^{\infty}\gamma(\lambda,xt)\varphi(t)\, dt,
\end{equation}
where $\varphi$ is a completely monotonic function which is integrable on $(0,\infty)$. Writing $\varphi=\mathcal L (\mu)$  it easily follows that
$$
\int_0^{\infty}\frac{d\mu(t)}{t}=\int_0^{\infty}\varphi(t)\, dt<\infty.
$$
Since
$$
(-1)^k\varphi^{(k)}(x)=\int_0^{\infty}t^ke^{-xt}\, d\mu(t), \quad k\geq 0,
$$
all derivatives of $\varphi$ are bounded on $(0,\infty)$ if and only if $\mu$ has moments of all orders. In this case we let $\varphi^{(k)}(0+)\equiv \lim_{t\to 0}\varphi^{(k)}(t)$.
\begin{prop}
\label{prop:asymp}
    Let $f$ be a bounded complete generalized Bernstein function of order $\lambda$ the form \eqref{eq:bounded}, where $\varphi$ and all its derivatives are bounded. Then $f$ has an asymptotic expansion
    $$
   f(x)=b+\int_0^{\infty}\varphi(t)\, dt-\sum_{k=0}^{n-2}\frac{(\lambda)_{k+1}}{(k+1)!}\frac{\varphi^{(k)}(0+)}{x^{k+1}}+(-1)^nQ_n(x),
    $$
    where the remainder $Q_n$ is a completely monotonic function of order $n-\lambda$. 
\end{prop}
{\it Proof.}
Since $\varphi$ is completely monotonic, bounded and integrable on $(0,\infty)$ we may define the completely monotonic function $\Phi$ as 
$$
\Phi(x)=\int_x^{\infty}\phi(t)\, dt, \quad x>0.
$$
A computation shows that, 
with $d\omega(t)=(1/t)d\mu(t)$,
$$
f(x)=b+
x^{\lambda}\int_0^{\infty}t^{\lambda-1}\Phi(t)e^{-xt}\, dt= b+
\Gamma(\lambda)x^{\lambda}\int_0^{\infty}\frac{d\omega (t)}{(x+t)^{\lambda}}.
$$
Since $\omega$ has moments of all orders there is an asymptotic expansion available (see \cite[Theorem 3.2 and Corollary 3.13]{KP2}). It reads
$$
x^{\lambda-1}\int_0^{\infty}\frac{d\omega (t)}{(x+t)^{\lambda}}=\sum_{k=0}^{n-1}\frac{(\lambda)_k(-1)^k}{k!}\frac{s_k(\omega)}{x^{k+1}}+(-1)^nR_n(x),
$$
where $s_k(\omega)$ is the $k$'th moment of $\omega$ and where $R_n$ is a completely monotonic function of order $n+1-\lambda$.
Therefore, 
$$
f(x)=b+\sum_{k=0}^{n-1}\frac{(\lambda)_k(-1)^k}{k!}\frac{s_k(\omega)}{x^{k}}+(-1)^nQ_n(x),
$$
is an asymptotic expansion of $f$ in which the remainder $Q_n$ is a completely monotonic function of order $n-\lambda$.

Since $\varphi^{(k)}(0+)=(-1)^ks_k(\mu)=(-1)^ks_{k+1}(\omega)$, the asymptotic expansion can be written in the form
$$
f(x)=b+\int_0^{\infty}\varphi(t)\, dt-\sum_{k=0}^{n-2}\frac{(\lambda)_{k+1}}{(k+1)!}\frac{\varphi^{(k)}(0+)}{x^{k+1}}+(-1)^nQ_n(x),
$$
and this completes the proof. \hfill $\square$

We conclude the paper with an example concerning the Lomax distribution.
\begin{ex}
    The distribution function $F_{\lambda}$ of a randomized Lomax distribution was considered in \cite[Example 3.11]{KP1} (see also the paragraph before \cite[Corollary 5.6]{KP4}). It can be written in the form
    $$
    F_{\lambda}(x)=\int_0^{\infty}\gamma(\lambda, xt)\frac{dt}{(1+t)^2},
    $$
    and thus Proposition \ref{prop:asymp} can be applied. For the completely monotonic function $\varphi(t)=(1+t)^{-2}$ it is easily seen that $\varphi^{(k)}(0)=(-1)^k(k+1)!$. Hence the asymptotic expansion is
    $$
    F_{\lambda}(x)=1+\sum_{k=1}^{n-1}(-1)^k\frac{(\lambda)_k}{x^k} +(-1)^nQ_n(x),
    $$
    where $Q_n$ is a completely monotonic function  of order $n-\lambda$.
\end{ex}

\noindent
Stamatis Koumandos\\
Department of Mathematics and Statistics\\
The University of Cyprus\\
P. O. Box 20537\\
1678 Nicosia, Cyprus\\
email: skoumand@ucy.ac.cy
\medskip

\vspace{0.1in}

\noindent
Henrik Laurberg Pedersen\\
Department of Mathematical Sciences\\
University of Copenhagen\\
Universitetsparken 5\\
DK-2100, Denmark\\
email: henrikp@math.ku.dk


\begin{thebibliography}{xx}

\bibitem{as} M.\ Abramowitz and I.A.\ Stegun: Handbook of
Mathematical Functions with formulas, Graphs and Mathematical
Tables. Dover, New York (1965)
\bibitem{aar}G.E.\ Andrews, R.\ Askey and R.\ Roy,
Special functions, Encyclopedia of Mathematics and its Applications, {\bf 71},
Cambridge University Press, Cambridge (1999).
\bibitem{BKP}C.\ Berg, S.\ Koumandos and H.L.\ Pedersen, Nielsen's beta function and some infinitely divisible
distributions, Math.\ Nach. {\bf 294} (2021) no.\ 3, 426--449. 
\bibitem{dlmf} NIST Digital Library of Mathematical Functions (dlmf.nist.gov).
 \bibitem{GR} I.S.~Gradshteyn and I.M.~Ryzhik, Table of Series and Products, Seventh Edition, Academic Press (2007).
\bibitem{ism} M.E.H.\ Ismail, L.\ Lorch and M.E.\ Muldoon, Completely
  monotonic functions associated with the Gamma function and its
  $q$-analogues. J. Math. Anal. Appl. \textbf{116},(1986), 1--9.

 \bibitem{kofte} D.A.\ Kofke, Comment on The incomplete beta function law for parallel tempering sampling of classical canonical systems (J.\ Chem.\ Phys.\ {\bf 120}, 4119 (2004)). J.\ Chem.\ Phys.\ {\bf 121} (2), pp.\ 1167.
\bibitem{KL} S.\ Koumandos and M.\ Lamprecht: Complete monotonicity and related properties of some special functions. Math. Comp. 82 (2013), no. 282, 1097--1120.
\bibitem{KP0}  S.\ Koumandos and H.L.\ Pedersen: Completely monotonic functions of positive order and asymptotic expansions of the logarithm of Barnes double gamma function and Euler’s gamma function. J.\ Math.\ Anal.\ Appl.\ {\bf 355} (2009), 33--40.
 \bibitem{KP1} S.\ Koumandos and H.L.\ Pedersen, Generalized Bernstein Functions, Math.\ Scand., {\bf 129} (2023), 93--116.
 \bibitem{KP2} S.\ Koumandos and H.L.\ Pedersen, On Asymptotic Expansions of Generalized Stieltjes Functions, Comput.\ Methods Funct.\ Theory, {\bf 15} (2015), 93--115.
 \bibitem{KP3} S.\ Koumandos and H.L.\ Pedersen,  On generalized Stieltjes Functions,   Constr. Approx., {\bf 50} (2019), no.\ 1, 129--144.
  \bibitem{KP4} S.\ Koumandos and H.L.\ Pedersen, Higher order Thorin-Bernstein functions, Results Math., {\bf 79} (2024), 22 pages. DOI: 10.1007/s00025-023-02040-z
\bibitem{kristiansen}G.K.~Kristiansen, A Proof of Steutel's Conjecture, \emph{Ann.\ Probab.} \textbf{22} (1994), 442--452.
 \bibitem{S}R.L.~Schilling, R.~Song and Z.~Vondracek, Bernstein Functions
Theory and Applications 2ND Ed. De Gruyter Studies in Mathematics, {\bf 37} (2012).

\end{thebibliography}
\end{document}